\documentclass[12pt,twoside]{article}

\usepackage{a4}
\usepackage{exscale}
\usepackage{fancyhdr}
\usepackage[centertags]{amsmath}
\usepackage{amsfonts}
\usepackage{amssymb}
\usepackage{amsthm}
\usepackage{ulem}
\usepackage{isolatin1}

\renewcommand{\sectionmark}[1]
                    {\markboth{Kapitel \thesection\ #1}{}}
\renewcommand{\sectionmark}[1]
                 {\markright{} }

\newtheorem{thm}{Theorem}
\newtheorem*{definition}{Definition}       
\newtheorem{lem}[thm]{Lemma}
\newtheorem{corollary}[thm]{Corollary}
\newtheorem{proposition}[thm]{Proposition}

\newtheorem{fremdersatz}{Theorem}

\newenvironment{satOrig}[1]
        {\pagebreak[2] \begin{fremdersatz} {\bf #1} \quad\sl}
        { \end{fremdersatz}}

\newenvironment{proofof}[1]
        {\pagebreak[2] \vspace{-1pt}{\bf Proof#1.}  }
        {\hfill $\blacksquare$ \vspace{2pt}}

\newenvironment{einr}{\parmod 
                      \begin{list}{}
                        {\setlength{\rightmargin}{0cm}
                         \setlength{\leftmargin}{0,75cm}
                         \setlength{\labelwidth}{0cm}
                         \setlength{\parsep}{1pt}
                         \setlength{\itemsep}{1pt}
                         \setlength{\topsep}{1pt}
                         \setlength{\partopsep}{0pt}
                         \setlength{\labelsep}{0cm}
                         \setlength{\listparindent}{0pt}
                         \setlength{\itemindent}{0pt}}
                      \item[] \ignorespaces}
                     {\unskip \end{list}}

\numberwithin{equation}{section}

\def\parmod{\parskip=2pt plus1pt minus1pt}

\def\ma{\mathcal{M}}

\def\l{\left}
\def\r{\right}

\newcommand{\fa}{\mathcal{F}}

\newcommand{\nat}{{\rm I\! N}}
\newcommand{\nato}{{\rm I\! N}_0}

\newcommand{\co}{{\mathbb C}}
\newcommand{\re}{{\mathbb R}}

\newcommand{\zet}{{\mathbb Z}}

\newcommand{\gl}{\left\{}
\newcommand{\gr}{\right\}}
\newcommand{\kl}{\left(}
\newcommand{\kr}{\right)}
\newcommand{\el}{\left[}
\newcommand{\er}{\right]}
\newcommand{\kj}{\overline}
\renewcommand{\l}{\left}
\renewcommand{\r}{\right}
\newcommand{\limn}{\lim_{n\to\infty}}
\newcommand{\ti}{\widetilde}
\newcommand{\In}{\subseteq}
\newcommand{\mi}{\setminus}
\newcommand{\abb}{\longrightarrow}
\newcommand{\convl}{\Longrightarrow}

\renewcommand{\rho}{\varrho}
\renewcommand{\phi}{\varphi}
\renewcommand{\epsilon}{\varepsilon}

\newcommand{\beq}{\begin{equation}}
\newcommand{\eeq}{\end{equation}}
\newcommand{\beqar}{\begin{eqnarray}}
\newcommand{\eeqar}{\end{eqnarray}}
\newcommand{\beqaro}{\begin{eqnarray*}}
\newcommand{\eeqaro}{\end{eqnarray*}}
\newcommand{\bsat}{\begin{theorem}}
\newcommand{\esat}{\end{theorem}}
\newcommand{\bsatorig}{\begin{satOrig}}
\newcommand{\esatorig}{\end{satOrig}}
\newcommand{\blem}{\begin{lemma}}
\newcommand{\elem}{\end{lemma}}
\newcommand{\bkor}{\begin{corollary}}
\newcommand{\ekor}{\end{corollary}}
\newcommand{\bbew}{\begin{proofof}}
\newcommand{\ebew}{\end{proofof}}

\def\bsat{\begin{thm}}
\def\esat{\end{thm}}
\def\blem{\begin{lem}}
\def\elem{\end{lem}}
\def\bkor{\begin{corollary}}
\def\ekor{\end{corollary}}
\def\bprop{\begin{proposition}}
\def\eprop{\end{proposition}}
\def\bdefin{\begin{definition}}
\def\edefin{\end{definition}}
\def\beinr{\begin{einr}}
\def\eeinr{\end{einr}}

\renewcommand{\rho}{\varrho}
\renewcommand{\phi}{\varphi}
\renewcommand{\epsilon}{\varepsilon}

\begin{document}

\thispagestyle{plain}

\begin{center}
{\LARGE\bf Quasi-normality induced by differential inequalities\\[18pt]}  

{\Large \it J\"urgen Grahl and Shahar Nevo}
\end{center}


\centerline{\bf Abstract} 

We show that the family $\fa_k$ of all meromorphic functions $f$ in a
domain $D$ satisfying 
$$\frac{|f^{(k)}|}{1+|f|}(z)\ge C \qquad \mbox{ for all }
z\in D$$ 
(where $k$ is a natural number and $C>0$) is quasi-normal. 
The proof relies mainly on the Zalcman-Pang rescaling method.  



\section{Introduction and statement of results}

According to Marty's theorem, a family $\fa$ of meromorphic functions
in a domain $D$ in the complex plane $\co$ is normal (in the sense of
Montel) if and only if the family $\gl f^\#\;:\;f\in\fa\gr$ of the
corresponding spherical derivatives $f^\#:=\frac{|f'|}{1+|f|^2}$ is
locally uniformly bounded in $D$.

In \cite{GrahlNevo-Spherical} we studied families of meromorphic
functions whose spherical derivatives are bounded away from zero and
proved the following counterpart to Marty's theorem. 

\bsatorig{}\label{MartyBelow}
Let $D\In\co$ be a domain and $C>0$. Let $\fa$ be a family
of functions meromorphic in $D$ such that 
$$f^\#(z)\ge C \qquad\mbox{ for all } z\in D \mbox{ and all }
f\in \fa.$$
Then $\fa$ is normal in $D$. 
\esatorig

Hence, the condition $\frac{|f'|}{1+|f|^2}(z)=f^\#(z)\ge C$ can be
considered as a differential inequality that induces normality. In
\cite{BarGrahlNevo}, \cite{ChenNevoPang}, \cite{GN-Marty},
\cite{GNP-NonExplicit} and \cite{LiuNevoPang} we studied more general
differential inequalities, involving higher derivatives, with respect
to the question whether they induce normality or at least
quasi-normality.

Before summarizing the main results from these studies, as far as they
are relevant in the context of the present paper, we would like to remind
the reader of the definition of quasi-normality and also to introduce
some notations.  

A family $\fa$ of meromorphic functions in a domain $D\In\co$ is said
to be {\it quasi-normal} if from each sequence $\gl f_n\gr_n$ in $\fa$
one can extract a subsequence which converges locally uniformly (with
respect to the spherical metric) on $D\mi E$ where the set $E$ (which
may depend on $\gl f_n\gr_n$) has no accumulation point in $D$. If the
exceptional set $E$ can always be chosen to have at most $q$ points,
yet for some sequence there actually occur $q$ such points, then we say
that $\fa$ is {\it quasi-normal of order} $q$.  

For $z_0\in\co$ and $r>0,$ we set
$\Delta(z_0,r):=\{z\in\co:|z-z_0|<r\} $.  By $\ma(D)$ we denote the
space of all functions meromorphic in a domain $D$. We write $P_f$ and
$Z_f$ for the set of poles resp.~for the set of zeros of a meromorphic
function $f$, and we use the notation ``$f_n\convl f$ (on $D$)'' to
indicate that the sequence $\{f_n\}_n$ converges to $f$ locally
uniformly in $D$ (with respect to the spherical metric or in the
holomorphic case with respect to the Euclidean metric). Finally, we
make use not only of the well-known Landau notation $a_n=O(b_n)$, but
also of $a_n=\Omega(b_n)$ which means that
$\liminf_{n\to\infty}\frac{|a_n|}{|b_n|}>0$. (The usage of $\Omega(.)$
is quite ambiguous in literature. Our definition is the one given by
D. Knuth.)

Now we can turn to the results on differential inequalities and
(quasi-)normality known so far. While \cite{GN-Marty} and
\cite{GNP-NonExplicit} dealt with generalizations of Marty's theorem
(more precisely with conditions of the form
$\frac{|f^{(k)}|}{1+|f|^\alpha}(z)\le C$), in \cite{BarGrahlNevo},
\cite{ChenNevoPang} and \cite{LiuNevoPang} the following extensions of
Theorem \ref{MartyBelow} were proved.

\bsatorig{} \label{PrevResults}
Let $k\ge 1$ and $j\ge 0$ be integers and $C>0$, $\alpha>1$
be real numbers. Let $\fa$ be a family of meromorphic functions in
some domain $D$. 
\begin{itemize}
\item[(a)] \cite{ChenNevoPang} 
If 
$$\frac{|f^{(k)}|}{1+|f|^\alpha}(z)\ge C \qquad \mbox{ for all }
z\in D \mbox{ and all } f\in\fa,$$
then $\fa$ is normal.
\item[(b)] \cite{LiuNevoPang} 
If 
$$\frac{|f'|}{1+|f|}(z)\ge C \qquad \mbox{ for all }
z\in D \mbox{ and all } f\in\fa,$$
then $\fa$ is quasi-normal.
\item[(c)] \cite{BarGrahlNevo}
If $k>j$ and
$$\frac{|f^{(k)}|}{1+|f^{(j)}|^\alpha}(z)\ge C \qquad \mbox{ for all }
z\in D \mbox{ and all } f\in\fa,$$
then $\fa$ is quasi-normal in $D$. If all functions in $\fa$ are
holomorphic, $\fa$ is quasi-normal of order at most $j-1$. (For $j=0$
and $j=1$ this means that it is normal.) 
\end{itemize}
\esatorig

In (a), one can expect only quasi-normality if $\alpha=1$, not
normality. This is demonstrated by the family $\gl z\mapsto nz^k\;: \;
n\in\nat\gr$ which is not normal at $z=0$.

As to (c), there is no reasonable way to extend it to the case
$\alpha=1$ if $j\ge 1$ as we have already pointed out in
\cite{BarGrahlNevo}: The sequence of the functions $f_n(z):=z^n-3^n$
satisfies $\frac{|f_n^{(k)}|}{1+|f_n^{(j)}|}\convl \infty$ on the
annulus $\gl z\in\co\;:\; 2<|z|<4\gr$ whenever $k>j\ge 1$, but $\gl
f_n\gr_n$ is not normal at any point $z$ with $|z|=3$, hence it isn't
quasi-normal. In fact, it isn't even $Q_\alpha$-normal for any ordinal
number $\alpha$. (For the exact definition of $Q_\alpha$-normality we
refer to \cite{Nevo-Transfinite}.) The same example also shows that
(a) cannot be extended to the case $\alpha<1$ (see
\cite[Theorem~2]{LiuNevoPang} for details).

So the question remains whether part (b) of Theorem \ref{PrevResults}
can be extended to higher derivatives. In this paper we answer this
question affirmatively. 

\bsat{}\label{mainresult}
Let $k$ be a natural number and $C>0$. Let $\fa_k$ be the family of
all meromorphic functions $f$ in some domain $D$ satisfying
$$\frac{|f^{(k)}|}{1+|f|}(z)\ge C \qquad \mbox{ for all } z\in D.$$
Then $\fa_k$ is quasi-normal in $D$. 
\esat

In the spirit of Bloch's heuristic principle, one might ask for a
corresponding result for functions meromorphic in the whole complex
plane. Indeed, such a result holds, but it is quite trivial: If
$f\in\ma(\co)$ would satisfy $\frac{|f^{(k)}|}{1+|f|}(z)\ge C>0$ for all
$z\in\co$, then $f^{(k)}$ would be a constant by Montel's theorem, so
$f$ would be a polynomial of degree $k$. This would imply
$\lim_{z\to\infty} \frac{|f^{(k)}|}{1+|f|}(z)=0$, a contradiction. So
there aren't any functions $f\in\ma(\co)$ at all which satisfy
$\frac{|f^{(k)}|}{1+|f|}(z)\ge C$ for all $z\in\co$.

In Theorem \ref{mainresult} the order of quasi-normality can be
arbitrarily large, even if we restrict to holomorphic functions, as
the sequence of the functions $f_n(z):= n(e^z-1)$ on the strip $D:=\gl
z\in\co: -1<{\rm Re}\,(z)<1\gr$ demonstrates. Obviously, every
subsequence of $\gl f_n\gr_n$ is not normal exactly at the points
$z_j=2j\pi i$ with $j\in\zet$, so $\gl f_n\gr_n$ is quasi-normal of
infinite order in $D$. But for every natural number $k$ we have
$$\frac{|f_n^{(k)}|}{1+|f_n|}(z)
=\frac{|e^z|}{\frac{1}{n}+|e^z-1|}\ge \frac{1}{e(2+e)}\qquad \mbox{
  for all } z\in D \mbox{ and all } n.$$

\section{Some lemmas}

The most essential tool in our proof is the famous rescaling lemma of
L.~Zalc\-man~\cite{zalcman} and X.-C.~Pang~\cite{pang89,pang90}. Here
we require the following version from \cite{pangzalc2000a}.

\blem[{\bf Zalcman-Pang Lemma}] \label{zalclemma} 
Let $\fa$ be a family of meromorphic functions in a domain $D$ all of
whose zeros have multiplicity at least $m$ and all of whose poles have
multiplicity at least $p$ and let $-p<\alpha<m$. If $\fa$ is not
normal at some $z_0\in D$, then there exist sequences $\gl
f_n\gr_n\In \fa$, $\gl z_n\gr_n \In D$ and $\gl \rho_n\gr_n \In (0,1)$
such that $\gl\rho_n\gr_n$ tends to 0, $\gl z_n\gr_n$ tends to $z_0$ and such
that the sequence $\gl g_n\gr_n$ defined by
$$g_n(\zeta):= \frac{1}{\rho_n^\alpha}\cdot f_n(z_n+\rho_n \zeta)$$
converges locally uniformly in $\co$ (with respect to the spherical
metric) to a non-constant function $g$ meromorphic in $\co$.
\elem

Furthermore, we need the following well-known normality criterion due
to Y. Gu \cite{Gu} which corresponds to Hayman's alternative. 

\blem\label{Gu} 
Let $k\ge 1$ be an integer. Then the family of all functions $f$
meromorphic in a domain $D\In \co$ satisfying $f(z)\ne 0$,
$f^{(k)}(z)\ne 1$ for every $z\in D$ is normal. 
\elem


In our proof the concept of  differential polynomials turns out to be
helpful to simplify notations. For the convenience of the reader we
first recall the notions relevant in this context: 

\bdefin{}
Let $D\In \co$ be an arbitrary domain. A mapping $M: \ma(D)\abb\ma(D)$ given by   
$$M[u]=a\cdot \prod_{\nu=1}^{d} u^{(k_\nu)} \quad \mbox{ for all }
u\in\ma(D)$$ 
with $d\in\nato$, $a\in\co\mi\gl0\gr$ and $k_1,\dots,k_d\in\nato$ is
called a {\bf differential monomial} of  {\bf degree} $\deg(M):=d$ and {\bf weight}
$w(M):=\sum_{\nu=1}^{d} (1+k_\nu)$. 


A sum $P:=M_1+\dots+M_p$ of differential monomials $M_1,\dots,M_p$ 
which are linearly independent over $\ma(D)$ 
is called a {\bf differential polynomial} of {\bf degree}
$\deg(P):=\max\gl \deg(M_1),\dots,\deg(M_p)\gr$ and {\bf weight}
$w(P):=\max\gl w(M_1),\dots,w(M_p)\gr$. Obviously, we have $\deg(P)\le
w(P)$. \label{degreeweight}

If $\deg(M_1)=\dots=\deg(M_p)=:d$, we call  $P$ {\bf homogeneous (of
  degree $d$)}. 

For every differential polynomial $P$ there exists a differential
polynomial $P'$ such that $P'[u](z)=(P[u])'(z)$ for all $u\in\ma(D)$ and all
$z\in D$. It is easy to see that $\deg(P')=\deg(P)$ and $w(P')\ge
w(P)$. 
\edefin

The following lemma is used in two slightly different situations
within our proof. 

\blem\label{DiffPolLem}
Let $m$ and $j$ be natural numbers. Then for any
meromorphic function $g$ in a domain $D$ we have
$$\kl\frac{g}{g^{(j)}}\kr^{(m)} =\frac{g^{(m)}}{g^{(j)}}
-\frac{g}{g^{(j)}}\cdot\frac{g^{(j+m)}}{g^{(j)}}
+\frac{g}{g^{(j)}}\cdot \frac{Q_m\el g^{(j)}\er}{(g^{(j)})^m}
+\sum_{\ell=1}^{m-1} \frac{g^{(m-\ell)}}{g^{(j)}}\cdot \frac{Q_\ell\el
  g^{(j)}\er}{(g^{(j)})^\ell}$$
where for $\ell=1,\dots,m$ the $Q_\ell$ are homogeneous differential
polynomials of degree $\ell$ (or possibly $Q_\ell\equiv0$), all
differential monomials in $Q_\ell$ have weight $2\ell$ and $Q_m[u]$
doesn't contain any higher derivatives than $u^{(m-1)}$.  
\elem

\bbew{}
First we show by induction that for any meromorphic function $f$ and
for all $\ell\ge 1$ 
\beq\label{DerivReciprocal}
\kl\frac{1}{f}\kr^{(\ell)}=-\frac{f^{(\ell)}}{f^2}+\frac{\ti{Q}_\ell[f]}{f^{\ell+1}}
\eeq
where $\ti{Q}_1\equiv0$ and for $\ell\ge 2$ the $\ti{Q}_\ell$ are
homogeneous differential polynomials of degree $\ell$, all of whose
differential monomials have weight $2\ell$ and do not contain any
higher derivatives than $u^{(\ell-1)}$. 

Indeed, this is obviously true for $\ell=1$, and if
(\ref{DerivReciprocal}) holds for some $\ell\ge1$, 
differentiating yields 
$$\kl\frac{1}{f}\kr^{(\ell+1)}=-\frac{f^{(\ell+1)}}{f^2}+
\frac{1}{f^{\ell+2}}\cdot\kl 2f^{(l)}f' f^{\ell-1}
+f\cdot\ti{Q}_\ell'[f]-(\ell+1)\cdot f'\cdot\ti{Q}_\ell[f]\kr.$$ 
So if we set 
$$\ti{Q}_{\ell+1}[f]:=2f^{(l)}f'
  f^{\ell-1}+f\cdot \ti{Q}_\ell'[f]-(\ell+1)\cdot f'\cdot\ti{Q}_\ell[f],$$
then (\ref{DerivReciprocal}) holds for $\ell+1$ instead of $\ell$, and 
$\ti{Q}_{\ell+1}$ has the desired properties. 

Now inserting
(\ref{DerivReciprocal}) into  
$$\kl\frac{g}{g^{(j)}}\kr^{(m)}=\sum_{\ell=0}^{m} {m\choose \ell} \cdot
g^{(m-\ell)} \cdot \kl\frac{1}{g^{(j)}}\kr^{(\ell)}$$
and noting that $\frac{g^{(j+\ell)}}{(g^{(j)})^2}$ can be written in the form
$$\frac{g^{(j+\ell)}}{(g^{(j)})^2}=\frac{(g^{(j)})^{\ell-1}\cdot g^{(j+\ell)}}{(g^{(j)})^{\ell+1}}$$ 
we immediately obtain the assertion of the lemma.  
\ebew

\section{Proof of Theorem \ref{mainresult}}
We apply induction. The quasi-normality of $\fa_1$ follows from
Theorem \ref{PrevResults} (b). 
\beinr
{\small 
In order to make our proof as self-contained as possible, we include a
simplified version of the proof of Theorem \ref{PrevResults} (b) here:
For every $f\in\fa_1$ we have $\l|\tfrac{f}{f'}(z)\r|\le \frac{1}{C}$
and $|f'(z)|\ge C$ for all $z\in D$. Thus $\gl
\tfrac{f}{f'}:f\in\fa_1\gr$ and $\gl f':f\in\fa_1\gr$ are normal in $D$
by Montel's theorem. So for a given sequence $\gl f_n\gr_n$ in $\fa_1$
we may assume (after moving to an appropriate subsequence which we
still denote by $\gl f_n\gr_n$) that $\gl\tfrac{f_n}{f'_n}\gr_n$
converges to a holomorphic limit function $h_0$ and that $\gl
f_n'\gr_n$ converges to some $d\in \ma(D)\cup\gl\infty\gr$.

If $h_0\not\equiv 0$, then we deduce that $\gl f_n\gr_n$ tends to $d
\cdot h_0$ locally uniformly in $D\mi Z_{h_0}$. (Here $d\cdot h_0$ is
understood to be $\equiv \infty$ if $d\equiv\infty$.) Since the zero
set $Z_{h_0}$ is isolated in $D$, $\gl f_n\gr_n$ turns out to be
quasi-normal in $D$.

Now we consider the case $h_0\equiv 0$ and show that in this case $\gl
f_n\gr_n$ is normal in the whole of $D$. We assume that  $\gl
f_n\gr_n$ is not normal at some point $z_0\in D$. Then from
Gu's Theorem (Lemma \ref{Gu}) and the fact that $f'_n$ omits the value
$\frac{C}{2}$ we can conclude that there is a sequence $\gl
z_n\gr_n$ tending to $z_0$ such that $f_n(z_n)=0$ for $n$ large
enough.
By Weierstra\ss' theorem we have 
$$0\equiv h_0\Longleftarrow \kl\frac{f_n}{f_n'}\kr'
=1-\frac{f_n f_n''}{(f_n')^2}$$
and in particular 
$$0=\limn \kl1-\frac{f_n f_n''}{(f_n')^2}\kr(z_n)=1,$$
a contradiction. This shows that $\gl f_n\gr_n$ is normal in $D$ for
$h_0\equiv 0$ and completes the proof of the theorem for $k=1$.  
}
\eeinr
Let some $k\ge 2$ be given and assume that it is already known that
(on arbitrary domains) each of the conditions 
$$\frac{|f'|}{1+|f|}(z)\ge C,\dots,
\frac{|f^{(k-1)}|}{1+|f|}(z)\ge C$$ 
implies quasi-normality.

Let $\gl f_n\gr_n$ be some sequence in $\fa_k$. As in the case $k=1$, from
$$|f_n^{(k)}(z)|\ge C \quad\mbox{ and } \quad
\frac{|f_n(z)|}{|f_n^{(k)}(z)|}\le\frac{1}{C} \qquad \mbox{ for all }
z\in D \mbox{ and all } n$$
and from Montel's theorem we conclude that
$\gl f_n^{(k)}\gr_n$ and $\gl \frac{f_n}{f_n^{(k)}}\gr_n$ are normal,
so after turning to an appropriate subsequence we may assume that $\gl
f_n^{(k)}\gr_n$ converges locally uniformly in $D$ to some limit
function $d\in \ma(D)\cup\gl\infty\gr$ and that $\gl
\frac{f_n}{f_n^{(k)}}\gr_n$ converges to some holomorphic limit
function $h_0$. Now we consider several cases.

{\bf Case 1:} There is an $m\in\gl 1,\dots,k-1\gr$ such that a
subsequence of $\gl\frac{f_n^{(m)}}{f_n^{(k)}}\gr_n$ converges to
$\infty$ locally uniformly in $D$. 

Without loss of generality we may assume that the sequence 
$\gl\frac{f_n^{(m)}}{f_n^{(k)}}\gr_n$ itself converges to $\infty$. Then
for any compact disk $\kj{\Delta(z_0,r)}\In D$ we have for $n$ large enough,
say for $n\ge n_0$
$$|f_n^{(m)}(z)|\ge |f_n^{(k)}(z)| \qquad\mbox{ for all } z\in
\kj{\Delta(z_0,r)}.$$
But this means 
$$\frac{|f_n^{(m)}|}{1+|f_n|}(z)\ge C \qquad\mbox{ for all } z\in
\kj{\Delta(z_0,r)} \mbox{ and all } n\ge n_0,$$
so $\gl f_n\gr_n$ is quasi-normal in $\Delta(z_0,r)$ by the induction
hypothesis. Since this holds at any point $z_0\in D$, we deduce the
quasi-normality of $\gl f_n\gr_n$ in the whole of $D$. 

{\bf Case 2:} $d\in \ma(D)$. 

Let $z_0$ be an arbitrary point in $D$ which is not a pole of
$d$. Then there is an $r>0$ and an $M<\infty$ such that
$\kj{\Delta(z_0,r)}\In D$ and $|d(z)|\le M-1$ for all $z\in
\kj{\Delta(z_0,r)}$. So for $n$ large enough we obtain 
$$\frac{M}{1+|f_n(z)|}\ge \frac{|f_n^{(k)}(z)|}{1+|f_n(z)|}\ge C,
\quad \mbox{ hence } \quad |f_n(z)|\le \frac{M}{C} \qquad \mbox{ for
  all } z\in \kj{\Delta(z_0,r)}.$$
Therefore $\gl f_n\gr_n$ is normal in
$\Delta(z_0,r)$ by Montel's theorem. This shows that $\gl f_n\gr_n$ is
normal in $D\mi P_d$. Since the poles of $d$ are isolated in $D$, we
obtain the quasi-normality of $\gl f_n\gr_n$ in $D$.

{\bf Case 3:} $h_0\not\equiv 0$. 

Here we can conclude that $\gl f_n\gr_n$ tends to $d\cdot h_0$ in $D\mi
Z_{h_0}$. Since $Z_{h_0}$ is isolated in $D$, again we have
quasi-normality of $\gl f_n\gr_n$.

{\bf Case 4:} $h_0\equiv0$, $d\equiv\infty$, and there doesn't exist
an $m\in\gl 1,\dots,k-1\gr$ such that a subsequence of
$\gl\frac{f_n^{(m)}}{f_n^{(k)}}\gr_n$ converges to $\infty$ locally
uniformly in $D$.

This case turns out to be the most recalcitrant one. To deal with it
we will extensively consider the functions
$\frac{f_n^{(\mu)}}{f_n^{(k)}}$ with $\mu<k$ and their behaviour at
certain points close to a given point of non-normality. It is useful
to note that these functions are holomorphic since $f_n^{(k)}$ has no
zeros and since poles of $f_n$ are zeros of
$\frac{f_n^{(\mu)}}{f_n^{(k)}}$.

{\bf Claim 1.} Let $\mu\in\gl 1,\dots,k-1\gr$ be given and assume that 
$$\frac{f_n^{(\mu+1)}}{f_n^{(k)}}\convl
h_{\mu+1},\dots,\frac{f_n^{(k-1)}}{f_n^{(k)}}\convl h_{k+1}$$
where the limit functions $h_{\mu+1},\dots,h_{k-1}$ are
holomorphic. (Of course, for $\mu=k-1$ this assumption is meant to be
void.) Then $\gl\frac{f_n^{(\mu)}}{f_n^{(k)}}\gr_n$ is normal.   

{\bf Proof of Claim 1.} Let's assume that
$\gl\frac{f_n^{(\mu)}}{f_n^{(k)}}\gr_n$ is not normal at some point
$z_0\in D$. Then, since the functions $\frac{f_n^{(\mu)}}{f_n^{(k)}}$
don't have any poles, we can apply the Zalcman-Pang lemma (Lemma
\ref{zalclemma}) with $\alpha=-\mu$. So (after turning to an appropriate
subsequence of $\gl f_n\gr_n$) we find sequences $\gl z_n\gr_n \In D$
and $\gl \rho_n\gr_n \In (0,1)$ such that $\gl\rho_n\gr_n$ tends to 0,
$\gl z_n\gr_n$ tends to $z_0$ and such that the sequence $\gl
g_n\gr_n$ defined by 
\beq\label{DefZPSeq}
g_n(\zeta):= \rho_n^\mu\cdot \frac{f_n^{(\mu)}}{f_n^{(k)}}(z_n+\rho_n
\zeta)
\eeq
converges locally uniformly in $\co$ to a non-constant entire function
$g$. We choose an arbitrary $\zeta_0\in\co$ such that $g(\zeta_0)\ne
0$ and consider the quantities
$$A_n:=\frac{f_n^{(\mu)}}{f_n^{(k)}}(z_n+\rho_n \zeta_0)$$
which satisfy
\beq\label{GrowthAn}
A_n\sim \frac{g(\zeta_0)}{\rho_n^\mu},
\eeq
in particular $\limn A_n=\infty$. 

{\bf Claim 1.1.}
\beq\label{Estim-NN1}
\frac{f_n^{(k+m)}}{f_n^{(k)}}(z_n+\rho_n \zeta_0)=O(\rho_n^{-m})
\qquad\mbox{ for } m=1,\dots,k-1.
\eeq
{\bf Proof of Claim 1.1.} 
By differentiating (\ref{DefZPSeq}) we obtain 
$$\rho_n^{\mu+1}\cdot
\kl\frac{f_n^{(\mu+1)}}{f_n^{(k)}}-\frac{f_n^{(\mu)}}{f_n^{(k)}}\cdot \frac{f_n^{(k+1)}}{f_n^{(k)}}\kr(z_n+\rho_n
\zeta)\convl g'(\zeta)\qquad\mbox{ in } \co.$$
Here $\gl \frac{f_n^{(\mu+1)}}{f_n^{(k)}}(z_n+\rho_n\zeta)\gr_n$ tends to 
the finite value $h_{\mu+1}(z_0)$ locally uniformly in $\co$. (Note
that this is also true for $\mu=k-1$ if we set $h_k:\equiv 1$.) So we 
also have 
$$\rho_n^{\mu+1}\cdot\frac{f_n^{(\mu)}}{f_n^{(k)}}\cdot \frac{f_n^{(k+1)}}{f_n^{(k)}}(z_n+\rho_n
\zeta)\convl -g'(\zeta) \qquad\mbox{ in } \co$$
and in particular
$$\limn \rho_n^{\mu+1}\cdot A_n\cdot \frac{f_n^{(k+1)}}{f_n^{(k)}}(z_n+\rho_n \zeta_0)=-g'(\zeta_0).$$
Combining this with (\ref{GrowthAn}) yields
$\frac{f_n^{(k+1)}}{f_n^{(k)}}(z_n+\rho_n \zeta_0)=O(\rho_n^{-1})$,
i.e. the assertion for $m=1$. 

Now assume that, for given $m\in\gl2,\dots,k-1\gr$,
\beq\label{Claim11-IndAssumpt}
\frac{f_n^{(k+r)}}{f_n^{(k)}}(z_n+\rho_n \zeta_0)=O(\rho_n^{-r})
\qquad\mbox{ holds for } r=1,\dots,m-1.
\eeq
If we differentiate (\ref{DefZPSeq}) $m$-times and apply Lemma
\ref{DiffPolLem} with $g:=f_n^{(\mu)}$ and $j=k-\mu$, condensing the
representation given there in a way sufficient for our current
purpose, we obtain  
\beqar\label{Claim1-Deriv}
g^{(m)}(\zeta)&\Longleftarrow& 
\rho_n^{m+\mu}\kl\frac{f_n^{(\mu)}}{f_n^{(k)}}\kr^{(m)}(z_n+\rho{_n\zeta})\nonumber\\ 
&=& \rho_n^{m+\mu} 
\kl-\frac{f_n^{(\mu)}}{f_n^{(k)}}\frac{f_n^{(k+m)}}{f_n^{(k)}}
+\sum_{\ell=0}^{m} \frac{f_n^{(\mu+m-\ell)}}{f_n^{(k)}}\cdot
\frac{Q_\ell\el f_n^{(k)}\er}{(f_n^{(k)})^\ell}\kr(z_n+\rho{_n\zeta}),
\eeqar
where $Q_0\equiv1$ and for $\ell\ge1$ the $Q_\ell$ are homogeneous
differential polynomials of degree $\ell$, all differential monomials
in $Q_\ell$ have weight $2\ell$ and $Q_m[u]$ doesn't contain any
higher derivatives than $u^{(m-1)}$.  

Here we can write 
$$\frac{Q_\ell\el f_n^{(k)}\er}{(f_n^{(k)})^\ell}
=\sum_{\sigma=1}^{s_\ell} c_{\ell,\sigma}\prod_{\nu=1}^{\ell}
\frac{f_n^{(d_{\ell,\sigma,\nu}+k)}}{f_n^{(k)}}$$
with certain constants $c_{\ell,\sigma}\in\co$, $s_\ell\ge0$ and
$d_{\ell,\sigma,\nu}\ge0$, where
$d_{\ell,\sigma,1}+\dots+d_{\ell,\sigma,\ell}=\ell$. For $\ell \le 
m-1$ all $d_{\ell,\sigma,\nu}$ are $\le m-1$, and in view of the
special property of $Q_m$ this also remains valid for $\ell=m$. So from
(\ref{Claim11-IndAssumpt}) we can conclude that
\beq\label{Claim11-DiffPol}
\frac{Q_\ell\el f_n^{(k)}\er}{(f_n^{(k)})^\ell}(z_n+\rho_n\zeta_0)=O(\rho_n^{-\ell})
\qquad \mbox{ for } \ell=0,\dots,m.
\eeq
Furthermore, from (\ref{Claim11-IndAssumpt}) and the assumption on the
existence of the limit functions $h_{\mu+1},\dots,h_{k-1}$ we can
conclude that for $\ell=0,\dots,m-1$ we have\footnote{Here, for
  applying (\ref{Claim11-IndAssumpt}) it is crucial that $m+\mu-\ell\le k+m-1$.} 
\beq\label{Claim11-LogDeriv}
\frac{f_n^{(\mu+m-\ell)}}{f_n^{(k)}}(z_n+\rho_n\zeta_0)=
\gl\begin{array}{ll} O(\rho_n^{-\mu-m+\ell+k}) & \mbox{ if } \mu+m-\ell>k,\\[5pt]
O(1) & \mbox{ if } \mu+m-\ell\le k.\end{array}\r.
\eeq
If we combine (\ref{Claim11-DiffPol}) and (\ref{Claim11-LogDeriv}), we
obtain for all $\ell=0,\dots,m-1$
$$\kl\frac{f_n^{(\mu+m-\ell)}}{f_n^{(k)}}\cdot 
\frac{Q_\ell\el f_n^{(k)}\er}{(f_n^{(k)})^\ell}\kr(z_n+\rho{_n\zeta_0})
=O(\rho_n^{-\mu-m+1}).$$
(Indeed, this is true in both cases $\mu+m-\ell>k$ and $\mu+m-\ell\le k$.)

Furthermore, in view of (\ref{GrowthAn}) and (\ref{Claim11-DiffPol})
we also have 
$$\kl\frac{f_n^{(\mu)}}{f_n^{(k)}}\cdot \frac{Q_m\el
  f_n^{(k)}\er}{(f_n^{(k)})^m}\kr(z_n+\rho{_n\zeta_0}) 
=O(\rho_n^{-\mu-m}).$$
Now we take a look again at (\ref{Claim1-Deriv}): Inserting
  $\zeta=\zeta_0$, we obtain for $n\to\infty$
\beqaro
&&\rho_n^{m+\mu} 
\kl\underbrace{A_n}_{=\Omega(\rho_n^{-\mu})}\cdot \frac{f_n^{(k+m)}}{f_n^{(k)}}
-\underbrace{\frac{f_n^{(\mu)}}{f_n^{(k)}}\cdot \frac{Q_m\el
  f_n^{(k)}\er}{(f_n^{(k)})^m}}_{=O(\rho_n^{-\mu-m})}
-\underbrace{\sum_{\ell=0}^{m-1} \frac{f_n^{(\mu+m-\ell)}}{f_n^{(k)}}\cdot
\frac{Q_\ell\el f_n^{(k)}\er}{(f_n^{(k)})^\ell}}_{=O(\rho_n^{-\mu-m+1})}\kr(z_n+\rho{_n\zeta_0})\\
&&\abb -g^{(m)}(\zeta_0),
\eeqaro
and we can conclude that $\frac{f_n^{(k+m)}}{f_n^{(k)}}(z_n+\rho_n
\zeta_0)=O(\rho_n^{-m})$. So by induction Claim 1.1 is proved for all
$m=1,\dots,k-1$. \hfill $\square$

{\bf Claim 1.2.} If there is some $m\in \gl 2,\dots,k-1\gr$ such that 
\beq\label{Claim12-Assump}
\frac{f_n^{(m)}}{f_n^{(k)}}(z_n+\rho_n \zeta_0)=\Omega(\rho_n^{-m}),
\eeq
then we can find a $j\in \gl 1,\dots,m-1\gr$ such that 
$$\frac{f_n^{(m-j)}}{f_n^{(k)}}(z_n+\rho_n\zeta_0)=\Omega(\rho_n^{-m+j}).$$

{\bf Proof of Claim 1.2.} Differentiating $\frac{f_n}{f_n^{(k)}}\convl 0$ $m$-times
and again applying Lemma \ref{DiffPolLem} (this time choosing a
slightly different way of condensing the representation given there
than in the proof of Claim 1.1), we obtain 
\beq\label{Claim12-Diff}
0\Longleftarrow\kl\frac{f_n}{f_n^{(k)}}\kr^{(m)}
=\frac{f_n^{(m)}}{f_n^{(k)}}
+\frac{f_n}{f_n^{(k)}}\cdot \frac{Q_m\el f_n^{(k)}\er}{(f_n^{(k)})^m}
+\sum_{\ell=1}^{m-1} \frac{f_n^{(m-\ell)}}{f_n^{(k)}}\cdot
\frac{Q_\ell\el f_n^{(k)}\er}{(f_n^{(k)})^\ell}
\eeq
where the $Q_\ell$ are homogeneous differential polynomials of degree $\ell$
and all differential monomials in $Q_\ell$ have weight $2\ell$. Again
we have
$$\frac{Q_\ell\el f_n^{(k)}\er}{(f_n^{(k)})^\ell}
=\sum_{\sigma=1}^{s_\ell} c_{\ell,\sigma}\prod_{\nu=1}^{\ell}
\frac{f_n^{(d_{\ell,\sigma,\nu}+k)}}{f_n^{(k)}}$$
with certain constants $c_{\ell,\sigma}\in\co$, $s_\ell\ge0$ and $d_{\ell,\sigma,\nu}\ge0$, where
$d_{\ell,\sigma,1}+\dots+d_{\ell,\sigma,\ell}=\ell\,(\le m\le
k-1)$. Therefore from Claim 1.1 we deduce 
$$\frac{Q_\ell\el f_n^{(k)}\er}{(f_n^{(k)})^\ell}(z_n+\rho_n\zeta_0)=O(\rho_n^{-\ell})
\qquad \mbox{ for } \ell=1,\dots,m.$$
Now inserting this estimate into (\ref{Claim12-Diff}) 
and using $\frac{f_n}{f_n^{(k)}}\convl0$ and (\ref{Claim12-Assump}) we obtain
$$\kl
\underbrace{\frac{f_n^{(m)}}{f_n^{(k)}}}_{=\Omega(\rho_n^{-m})}
+\underbrace{\frac{f_n}{f_n^{(k)}}}_{\abb 0}\cdot 
\underbrace{\frac{Q_m\el f_n^{(k)}\er}{(f_n^{(k)})^m}}_{=O(\rho_n^{-m})}
+\sum_{\ell=1}^{m-1} \frac{f_n^{(m-\ell)}}{f_n^{(k)}}\cdot
\underbrace{\frac{Q_\ell\el f_n^{(k)}\er}{(f_n^{(k)})^\ell}}_{=O(\rho_n^{-\ell})}\kr(z_n+\rho_n\zeta_0)
\abb0 \quad (n\to\infty)$$
which shows us that there must be a $j\in \gl 1,\dots,m-1\gr$ such that 
$$\frac{f_n^{(m-j)}}{f_n^{(k)}}(z_n+\rho_n\zeta_0)=\Omega(\rho_n^{-m+j}).$$
This proofs Claim 1.2. \hfill $\square$

Now we can complete the proof of Claim 1. 

If $\mu\ge 2$, then from Claim 1.2 and the fact that (\ref{Claim12-Assump})
holds for $m=\mu$ by our choice of $\zeta_0$ (see (\ref{GrowthAn})) we
deduce that for a suitable $j\in \gl 1,\dots,\mu-1\gr$ we have
$$\frac{f_n^{(\mu-j)}}{f_n^{(k)}}(z_n+\rho_n
\zeta_0)=\Omega(\rho_n^{-\mu+j}).$$
Iterated application of Claim 1.2 (with $m=\mu-j$ and so on) shows that
$$\frac{f_n'}{f_n^{(k)}}(z_n+\rho_n \zeta_0)=\Omega(\rho_n^{-1}),$$
and this of course remains valid also for $\mu=1$ (again by our choice
of $\zeta_0$). 

But on the other hand, from $h_0\equiv 0$ by Weierstra\ss{}' theorem
we conclude that 
$$\frac{f_n'}{f_n^{(k)}}-\frac{f_n}{f_n^{(k)}}\cdot
\frac{f_n^{(k+1)}}{f_n^{(k)}}\convl 0,$$
hence, using Claim 1.1
$$\underbrace{\frac{f_n'}{f_n^{(k)}}(z_n+\rho_n\zeta_0)}_{=\Omega(\rho_n^{-1})}
-\underbrace{\frac{f_n}{f_n^{(k)}}(z_n+\rho_n\zeta_0)}_{\abb 0}\cdot
\underbrace{\frac{f_n^{(k+1)}}{f_n^{(k)}}(z_n+\rho_n\zeta_0)}_{=O(\rho_n^{-1})}\abb
0\qquad (n\to\infty).$$
This yields a contradiction, thus completing the proof of Claim~1. \hfill $\blacksquare$ 

Applying Claim 1 with $m=k-1$ (where the additional assumption is
void), we first see that $\gl\frac{f_n^{(k-1)}}{f_n^{(k)}}\gr_n$ is
normal, and without loss of generality we may assume that it tends to
some limit function $h_{k-1}$. By the general assumption in Case 4 we have
$h_{k-1}\not\equiv \infty$. So we can apply Claim 1 once more, with
$m=k-2$, to deduce the normality of
$\gl\frac{f_n^{(k-2)}}{f_n^{(k)}}\gr_n$. Continuing in this way, we
recursively find that (after turning to suitable subsequences) all sequences 
$\gl\frac{f_n^{(m)}}{f_n^{(k)}}\gr_n$ with $m=1,\dots,k-1$ tend to
holomorphic limit functions $h_m$. Furthermore, we know that
$\gl\frac{f_n}{f_n^{(k)}}\gr_n$ tends to $h_0\equiv0$. So by Weierstra\ss{}'
theorem we also have 
\beq\label{Weierstr}
\frac{f_n^{(m+1)}}{f_n^{(k)}}-\frac{f_n^{(m)}}{f_n^{(k)}}\cdot \frac{f_n^{(k+1)}}{f_n^{(k)}}
=\kl\frac{f_n^{(m)}}{f_n^{(k)}}\kr'\convl h_m' \quad\mbox{ for }
m=0,\dots,k-1.
\eeq
Assume that $h_{m_0}\not\equiv 0$ for some $m_0\in
\gl1,\dots,k-1\gr$. Then from (\ref{Weierstr}) we deduce that $\gl
\frac{f_n^{(k+1)}}{f_n^{(k)}}\gr_n$ tends to some meromorphic limit
function $L$ in $D\mi Z_{h_{m_0}}$. Again by (\ref{Weierstr}) we
obtain in $D\mi Z_{h_{m_0}}$
$$h_{m+1}-h_m\cdot L=h_m' \quad\mbox{ for } m=0,\dots,k-1.$$
From this and from $h_0\equiv 0$ we successively conclude that
$h_1\equiv0$, $h_2\equiv0,\dots,h_{m_0}\equiv 0$, a contradiction.  

So all limit functions $h_0,\dots,h_{k-1}$ are 0, i.e.~it only remains
to consider the situation that (after the usual extraction of
appropriate subsequences) each of the sequences
$\gl\frac{f_n}{f_n^{(k)}}\gr_n,\gl\frac{f_n'}{f_n^{(k)}}\gr_n,\dots,
\gl\frac{f_n^{(k-1)}}{f_n^{(k)}}\gr_n$ tends to 0 locally uniformly in
$D$.  Here, by Weierstra\ss{}' theorem we obtain for $m=1,\dots,k$
\beqaro
0&\Longleftarrow&\kl\frac{f_n^{(k-m)}}{f_n^{(k)}}\kr^{(m)}
=\sum_{j=0}^{m} {m\choose j} \cdot
f_n^{(k-j)}\cdot\kl\frac{1}{f_n^{(k)}}\kr^{(j)}
=\sum_{j=0}^{k} {m\choose j} \cdot f_n^{(k-j)}\cdot\kl\frac{1}{f_n^{(k)}}\kr^{(j)}.
\eeqaro
(Note that ${m\choose j}=0$ for $j>m$.)

Now we use the fact that the truncated (non-symmetric) Pascal matrix
$$A_k:=\Biggl({m\choose j}\Biggr)_{m=1,\dots,k\quad\atop
  j=0,\dots,k-1}
=\begin{pmatrix}
1 & 1 & 0 & 0 & \dots & 0 & 0 \\
1 & 2 & 1 & 0 & \dots & 0 & 0 \\
1 & 3 & 3 & 1 & \dots & 0 & 0 \\
\vdots & \vdots & \vdots & \vdots & & \vdots & \vdots \\
1 & k-1 & {k-1\choose 2} & {k-1\choose 3} & \dots & {k-1\choose k-2} &  1 \\[5pt]
1 & k & {k\choose 2} & {k\choose 3} & \dots & {k\choose k-2} &
  {k\choose k-1} 
\end{pmatrix}$$
is regular. This is a special case of a result of S.~Kersey
(\cite[Theorem 1.1]{Kersey}) which states that such a submatrix of the
Pascal matrix is regular if and only if there aren't any zeros on the
diagonal\footnote{The regularity can also be proved by induction, by showing that
  $A_k$ can be transformed by Gaussian elimination into the $k\times
  k$-matrix  
$$\begin{pmatrix}
1 & 1 & 0 & 0 & \dots & 0 & 0 \\
0 & 1 & 1 & 0 & \dots & 0& 0 \\
0 & 0 & 1 & 1 & \dots & 0& 0 \\
\vdots & \vdots & \vdots & \vdots & & \vdots & \vdots \\
0 & 0 & 0 & 0 & \dots & 1 & 1 \\
0 & 0 & 0 & 0 & \dots & 0 & 1 
\end{pmatrix}$$
which is obviously regular. We omit the details.}. 

Therefore there exists a vector $v\in\re^k$ such that $v^T
A_k=(1,0,\dots,0)$, i.e. 
$$\sum_{m=1}^{k} v_m\cdot {m\choose j}=
\gl\begin{array}{ll}
1 & \mbox{ for } j=0,\\
0 & \mbox{ for } j=1,\dots,k-1.
\end{array}\r.$$
This yields
\beqar\label{PascalMatrixApplied}
0&\Longleftarrow&
\sum_{m=1}^{k} v_m\cdot \sum_{j=0}^{k} {m\choose j}
\cdot f_n^{(k-j)}\cdot\kl\frac{1}{f_n^{(k)}}\kr^{(j)}\nonumber\\
&=&\sum_{j=0}^{k-1} \sum_{m=1}^{k} v_m\cdot  {m\choose j}
\cdot f_n^{(k-j)}\cdot\kl\frac{1}{f_n^{(k)}}\kr^{(j)}
+v_k\cdot f_n\cdot\kl\frac{1}{f_n^{(k)}}\kr^{(k)}\nonumber\\
&=& 1+v_k\cdot f_n\cdot\kl\frac{1}{f_n^{(k)}}\kr^{(k)}.
\eeqar
If $\gl f_n\gr_n$ would be not normal at some point $z_0\in D$, then
from Gu's theorem (Lemma \ref{Gu}) and the fact that the $f_n^{(k)}$ omit
the value $\tfrac{C}{2}\ne0$ we could conclude that there is a
sequence $\gl z_n\gr_n$ tending to $z_0$ such that $f_n(z_n)=0$ for
$n$ large enough. Inserting this into (\ref{PascalMatrixApplied})
yields 
$$0= \limn\kl 1+v_k\cdot f_n\cdot\kl\frac{1}{f_n^{(k)}}\kr^{(k)}\kr(z_n)=1,$$
a contradiction. 

This completes the proof of Theorem \ref{mainresult}. 






\vspace{10pt}
\parbox{80mm}{\it J\"urgen Grahl\\
University of W\"urzburg \\
Department of Mathematics  \\     
97074 W\"urzburg\\
Germany\\
e-mail: grahl@mathematik.uni-wuerzburg.de}
\hfill\parbox{75mm}{\it Shahar Nevo \\
Bar-Ilan University\\
Department of Mathematics\\
Ramat-Gan 52900\\
Israel\\
e-mail: nevosh@math.biu.ac.il}

\end{document}